\documentclass[12pt]{amsart}
\usepackage{graphicx}
\usepackage{amssymb}
\usepackage{amsfonts}
\usepackage{amsmath}
\usepackage{array}
\usepackage{rotating}
\usepackage{epsfig}

\newtheorem{example}{Example}[section]

\newtheorem{theorem}[example]{Theorem}

\newtheorem{corollary}[example]{Corollary}

\newtheorem{proposition}[example]{Proposition}

\def\Proof{\noindent \it Proof -- \rm}
\def\qed{\hspace{3.5mm} \hfill \vbox{\hrule height 3pt depth 2 pt width 2mm}
\bigskip}

\def\QSym{{\it QSym}}          
\def\NCSF{{\bf Sym}}           
\def\FQSym{{\bf FQSym}}        
\def\FSym{{\bf FSym}}          
\def\PBT{{\bf PBT}}            
\def\CQSym{{\bf CQSym}}        
\def\CPQSym{{\it PQSym}}       
\def\EQSym{{\it EQSym}}        
\def\ESym{{\bf ESym}}          
\def\PiQSym{{\it {\Pi}QSym}}   
\def\CCQSym{{\it CQSym}}       
\def\SGQSym{{\it \SG QSym}}    
\def\SGSym{{\bf \SG Sym}}      
\def\PhiSym{{\bf {\Phi}Sym}}   

\def\MEF{{\it M}}              
\def\SEF{{\bf S}}              
\def\uq{{\it U}}               
\def\upi{{\it U}}              
\def\ul{{\it u}}               
\def\Mper{{\it M}}             
\def\Sper{{\bf S}}             
\def\Mpa{{\it M}}              

\def\csupp{{\rm supp}}         
\def\picyc{{\shuffle}}         
\def\picycgen{{\smile}}        

\def\ncbinomial#1#2{\left[\,\begin{matrix}#1 \cr #2\end{matrix}\,\right]}
\def\myphi{{\bf {\bf\phi}}} 

\def\inv{{\rm inv}}     

\def\sconc{\bullet}     
\def\Std{{\rm Std}}     
\def\cstd{{\rm cstd}}   

\def\<{\langle}
\def\>{\rangle}
\def\ZZ{{\mathbb Z}}    
\def\CC{{\mathbb C}}    
\def\KK{{\mathbb K}\, } 

\def\park{{\bf a}} 
\def\parkc{{\pi}} 
\def\tr{\operatorname{tr}}

\def\F{{\bf F}}         
\def\S{{\bf S}}         
\def\SG{{\mathfrak S}}  
\def\gr{{\rm gr}}

\def\tensor{\otimes}

\def\End{\operatorname{End}} 

\def\PF{{\rm PF}}   

\def\shuff#1#2{\mathbin{
\hbox{\vbox{ \hbox{\vrule \hskip#2 \vrule height#1 width 0pt
}%
\hrule}%
\vbox{ \hbox{\vrule \hskip#2 \vrule height#1 width 0pt
\vrule }%
\hrule}%
}}}

\def\qbin#1#2{\begin{bmatrix} #1 \\ #2\end{bmatrix}}

\def\shuf{{\mathchoice{\shuff{7pt}{3.5pt}}%
{\shuff{6pt}{3pt}}%
{\shuff{4pt}{2pt}}%
{\shuff{3pt}{1.5pt}}}}%
\def\shuffle{\,\shuf\,}

\title[Commutative Hopf algebras of permutations and trees]%
{Commutative Hopf algebras\\ of permutations and trees}

\author[F. Hivert, J.-C.~Novelli, and J.-Y.~Thibon]%
{Florent Hivert, Jean-Christophe Novelli, and Jean-Yves Thibon}

\address[] {Institut Gaspard Monge, Universit\'e de Marne-la-Vall\'ee \\
5 Boulevard Descartes \\Champs-sur-Marne \\77454 Marne-la-Vall\'ee cedex 2 \\
FRANCE}
\email[Florent Hivert]{hivert@univ-mlv.fr}
\email[Jean-Christophe Novelli]{novelli@univ-mlv.fr}
\email[Jean-Yves Thibon]{jyt@univ-mlv.fr} 
\date{}

\begin{document}

\begin{abstract}
We propose several constructions of commutative or cocommutative Hopf algebras
based on various combinatorial structures, and investigate the relations
between them. A commutative Hopf algebra of permutations is obtained by a
general construction based on graphs, and its non-commutative dual is realized
in three different ways, in particular as the Grossman-Larson algebra of heap
ordered trees.
Extensions to endofunctions, parking functions, set partitions, planar binary
trees and rooted forests are discussed. Finally, we introduce one-parameter
families interpolating between different structures constructed on the
same combinatorial objects.
\end{abstract}

\maketitle

\section{Introduction}

Many examples of Hopf algebras based on combinatorial structures are known.
Among these, algebras based on permutations and planar binary trees play a
prominent role, and arise in seemingly unrelated
contexts~\cite{MR,LR1,NCSF6,BF}.
As Hopf algebras, both are noncommutative and non cocommutative, and in fact
self-dual.

More recently, cocommutative Hopf algebras of binary trees and permutations
have been constructed~\cite{NT1,AS}. In~\cite{NT1}, binary trees arise as sums
over rearrangements classes in an algebra of parking functions, while
in~\cite{AS}, cocommutative Hopf algebras are obtained as the graded
coalgebras associated with coradical fitrations.

In~\cite{NTT}, a general method for constructing commutative Hopf algebras
based on various kind of graphs has been presented.  The aim of this note is
to investigate Hopf algebras based on permutations and trees constructed by
the method developed in~\cite{NTT}. These commutative algebras are, by
definition, realized in terms of polynomials in an infinite set of doubly
indexed indeterminates. The dual Hopf algebras are then realized by means of
non commutative polynomials in variables $a_{ij}$. We show that these first
resulting algebras are isomorphic (in a non trivial way) to the duals of those
of~\cite{AS}, and study some generalizations such as endofunctions, parking
functions, set partitions, trees, forests, and so on.

The possibility to obtain in an almost systematic way commutative, and in
general non cocommutative, versions of the usual combinatorial Hopf algebras
leads us to conjecture that these standard versions should be considered as
some kind of quantum groups, \emph{i.e.}, can be incorporated into
one-parameter families, containing an enveloping algebra and its dual for
special values of the parameter.
A few results supporting this point of view are presented in the final
section.

\medskip
In all the paper, $\KK$ will denote a field of characteristic zero.
All the notations used here is as in~\cite{NCSF1, NTT}.

\section{A commutative Hopf algebra of endofunctions}
\label{eqsym}

Permutations can be regarded in an obvious way as labelled and oriented graphs
whose connected components are cycles.
Actually, arbitrary \emph{endofunctions} (functions from $[n]:=\{1,\ldots,n\}$
to itself) can be regarded as labelled graphs, connecting $i$ with $f(i)$ for
all $i$ so as to fit in the framework of~\cite{NTT}, where a general process
for building Hopf algebras of graphs is described.

In the sequel, we identify an endofunction $f$ of $[n]$ with the word
\begin{equation}
w_f=f(1)f(2)\cdots f(n) \in [n]^n.
\end{equation}

Let $\{x_{i\,j} \,|\, i,j\geq1\}$ be an infinite set of commuting
indeterminates, and let $\mathcal J$ be the ideal of
$R=\KK[x_{i\,j} \,|\, i,j\geq1]$ generated by the relations
\begin{equation}
x_{i\,j} x_{i\,k}=0
\quad \text{for all $i,j,k$.}
\end{equation}
For an endofunction $f:[n]\to [n]$, define
\begin{equation}
\label{mf-eqs}
\MEF_f := \sum_{i_1 < \cdots < i_n}
            x_{i_1\, i_{f(1)}}\cdots x_{i_n\, i_{f(n)}},
\end{equation}
in $R/{\mathcal J}$.

It follows from~\cite{NTT}, Section~4, that

\begin{theorem}
The $\MEF_f$ span a subalgebra $\EQSym$ of the commutative algebra
$R/{\mathcal J}$.
More precisely, there exist non-negative integers $C_{f,g}^{h}$
such that
\begin{equation}
\MEF_f \MEF_g = \sum_{h} C_{f,g}^{h} \MEF_h.
\end{equation}
\qed
\end{theorem}

\begin{example}
{\rm
\begin{equation}
\MEF_{1} \MEF_{22} = \MEF_{133} + \MEF_{323} + \MEF_{223}.
\end{equation}
\begin{equation}
\MEF_{1} \MEF_{331} = \MEF_{1442} + \MEF_{4241} + \MEF_{4431} + \MEF_{3314}.
\end{equation}
\begin{equation}
\MEF_{12} \MEF_{21} = \MEF_{1243} + \MEF_{1432} + \MEF_{4231} + \MEF_{1324} +
\MEF_{3214} + \MEF_{2134}.
\end{equation}
\begin{equation}
\label{mf12-22}
\MEF_{12} \MEF_{22} = \MEF_{1244} + \MEF_{1434} + \MEF_{4234} + \MEF_{1334} +
\MEF_{3234} + \MEF_{2234}.
\end{equation}
}
\end{example}

The \emph{shifted concatenation} of two endofunctions $f:[n]\to[n]$ and
$g:[m]\to[m]$ is the endofunction
$h:=f\bullet g$ of $[n+m]$ such that $w_h :=w_f\bullet w_g$, that is
\begin{equation}
\left\{
\begin{aligned}
h(i)=f(i)     \ & \text{if} &\ i\leq n\\ 
h(i)=n+g(i-n) \ & \text{if} &\ i > n\\ 
\end{aligned}
\right.
\end{equation}

We can now give a combinatorial interpretation of the coefficient
$C_{f,g}^{h}$: if $f:[n]\to[n]$ and $g:[m]\to[m]$, this coefficient is the
number of permutations $\tau$ in the shuffle product
$(1\ldots p)\shuffle (p+1\ldots p+n)$ such that
\begin{equation}
h = 
\tau^{-1} \circ (f\bullet g) \circ \tau.
\end{equation}
For example, with $f=12$ and $g=22$, one finds the set (see
Equation~(\ref{mf12-22}))
\begin{equation}
\label{mf12-22a}
\{1244, 1434, 4234, 1334, 3234, 2234 \}.
\end{equation}

Now, still following~\cite{NTT}, define a coproduct by

\begin{equation}
\label{coprodM-endof}
\Delta \MEF_h := \sum_{(f,g) ; f\sconc g=h} 
\MEF_{f} \tensor \MEF_{g}.
\end{equation}
This endows $\EQSym$ with a (commutative, non cocommutative) Hopf algebra
structure.

\begin{example}
\begin{equation}
\Delta \MEF_{626124}  = \MEF_{626124} \tensor 1 + 1\tensor \MEF_{626124}.
\end{equation}
\begin{equation}
\Delta \MEF_{4232277} = \MEF_{4232277} \tensor 1 +
                        \MEF_{42322}\tensor \MEF_{22} +
                        1\tensor \MEF_{4232277}.
\end{equation}
\end{example}

Define a \emph{connected} endofunction as a function that cannot be obtained
by non trivial shifted concatenation. Then, the definition of the coproduct
of the $\MEF_f$ implies

\begin{proposition}
If $(\SEF^f)$ denotes the dual basis of $(\MEF_f)$, the graded
dual $\ESym := \EQSym^*$ is free over the set
\begin{equation}
\{ \SEF^f \,|\, f \text{\ connected} \}.
\end{equation}
\end{proposition}

Indeed, Equation~(\ref{coprodM-endof}) is equivalent to
\begin{equation}
\SEF^f \SEF^g = \SEF^{f\sconc g}.
\end{equation}
\qed

Now, $\ESym$ being a graded connected cocommutative Hopf algebra, it follows
from the Milnor-Moore theorem that
\begin{equation}
\ESym = U(L),
\end{equation}
where $L$ is the Lie algebra of its primitive elements.
Let us now prove
\begin{theorem}
As a graded Lie algebra, the primitive Lie algebra $L$ of $\ESym$ is free over
a set indexed by connected endofunctions.
\end{theorem}

\Proof
Assume it is the case. By standard arguments on generating series, one finds
that the number of generators of $L$ in degree $n$ is equal to the number of
algebraic generators of $\ESym$ in degree $n$, parametrized for example by
connected endofunctions (series beginning by
$(1,3,20,197,2511,38924,\ldots)$).
We will now show that $L$ has at least this number of generators and that
those generators are algebraically independent, determining completely the
dimensions of the homogeneous components $L_n$ of $L$ (series beginning by
$(1,3,23,223,2800,42576,\ldots)$).
Following Reutenauer~\cite{Reu} p.~58, denote by $\pi_1$ the Eulerian
idempotent, that is, the endomorphism of $\ESym$ defined by
$\pi_1=\log^*(Id)$. It is obvious, thanks to the definition of $\S^f$ that
\begin{equation}
\pi_1(\S^{f}) = \S^f + \cdots,
\end{equation}
where the dots stand for terms $\SEF^g$ such that $g$ is not connected.
Since the $\SEF^f$ associated with connected endofunctions are
independent, the dimension of $L_n$ is at least equal to the
number of connected endofunctions of size $n$. So $L$ is free over a set of
primitive elements parametrized by connected endofunctions.
\qed


There are many Hopf subalgebras of $\EQSym$ which can be defined by imposing
natural restrictions to maps: being bijective (see Section~\ref{sgqsym}),
idempotent ($f^2=f$), involutive ($f^2=id$), or more generally the Burnside
classes ($f^p=f^q$), and so on.
We shall start with the Hopf algebra of permutations.

\section{A commutative Hopf algebra of permutations}
\label{sgqsym}

\subsection{The Hopf algebra of bijective endofunctions}

Let us define $\SGQSym$ as the subalgebra of $\EQSym$ spanned by the
\begin{equation}
\Mper_\sigma = \sum_{i_1 < \cdots < i_n}
            x_{i_1\, i_{\sigma(1)}}\cdots x_{i_n\, i_{\sigma(n)}},
\end{equation}
where $\sigma$ runs over bijective endofunctions, \emph{i.e.}, permutations.
Note that $\SGQSym$ is also isomorphic to the image of $\EQSym$ in the
quotient of $R/{\mathcal J}$ by the relations

\begin{equation}
x_{i\,k}x_{j\,k}=0 \quad \text{for all $i,j,k$.}
\end{equation}

By the usual argument, it follows that

\begin{proposition}
The $\Mper_\sigma$ span a Hopf subalgebra $\SGQSym$ of the commutative Hopf
algebra $\EQSym$.
\qed
\end{proposition}

As already mentionned, there exist non-negative integers
$C_{\alpha,\beta}^{\gamma}$ such that
\begin{equation}
\Mper_\alpha \Mper_\beta = \sum_{\gamma} C_{\alpha,\beta}^{\gamma}
\Mper_\gamma.
\end{equation}

The combinatorial interpretation of the coefficients $C_{f,g}^h$ seen in
Section~\ref{eqsym} can be reformulated in the special case of permutations.
Write $\alpha$ and $\beta$ as a union of disjoint cycles. Split the set
$[n+m]$ into a set $A$ of $n$ elements, and its complement $B$, in
all possible ways. For each splitting, apply to $\alpha$ (resp. $\beta$) in
$A$ (resp. $B$) the unique increasing morphism of alphabets from
$[n]$ to $A$ (resp. from $[m]$ to $B$) and return the list of permutations
with the resulting cycles.
On the example $\alpha=12$ and $\beta=321$, this yields
\begin{equation}
\label{12-321b}
\begin{split}
(1)(2)(53)(4), (1)(3)(52)(4), (1)(4)(52)(3), (1)(5)(42)(3), (2)(3)(51)(4),\\
(2)(4)(51)(3), (2)(5)(41)(3), (3)(4)(51)(2), (3)(5)(41)(2), (4)(5)(31)(2).
\end{split}
\end{equation}
This set corresponds to the permutations and multiplicities of
Equation~(\ref{12-321c}).

A third interpretation of this product comes from the dual coproduct point of
view: $C_{\alpha,\beta}^{\gamma}$ is the number of ways of getting
$(\alpha,\beta)$ as the standardized words of pairs $(a,b)$ of two
complementary subsets of cycles of $\gamma$.
For example, with $\alpha=12$, $\beta=321$, and $\gamma=52341$, one has three
solutions for the pair $(a,b)$, namely
\begin{equation}
((2)(3), (4)(51)), ((2)(4), (3)(51)) , ((3)(4), (2)(51)),
\end{equation}
which is coherent with Equations~(\ref{12-321b}) and (\ref{12-321c}).

\begin{example}
\begin{equation}
\Mper_{12\cdots n} \Mper_{12\cdots p} =
  \binom{n+p}{n} \Mper_{12\cdots (n+p)}.
\end{equation}
\begin{equation}
\Mper_{1} \Mper_{21} = \Mper_{132} + \Mper_{213} + \Mper_{321}.
\end{equation}
\begin{equation}
\Mper_{12} \Mper_{21} = \Mper_{1243} + \Mper_{1324} + \Mper_{1432} +
\Mper_{2134} + \Mper_{3214} + \Mper_{4231}.
\end{equation}
\begin{equation}
\label{12-321c}
\Mper_{12} \Mper_{321} = \Mper_{12543} + \Mper_{14325} + 2 \Mper_{15342} +
\Mper_{32145} + 2 \Mper_{42315} + 3 \Mper_{52341}.
\end{equation}
\begin{equation}
\begin{split}
\Mper_{21} \Mper_{123} &=\
      \Mper_{12354} + \Mper_{12435} + \Mper_{12543} + \Mper_{13245} +
      \Mper_{14325}\\
  &+\ \Mper_{15342} + \Mper_{21345} + \Mper_{32145} + \Mper_{42315} +
      \Mper_{52341}.
\end{split}
\end{equation}
\begin{equation}
\begin{split}
\Mper_{21} \Mper_{231} &=\
   \Mper_{21453} + \Mper_{23154} + \Mper_{24513} + \Mper_{25431} +
   \Mper_{34152}\\
 &+\ \Mper_{34521} + \Mper_{35412} + \Mper_{43251} + 
   \Mper_{43512} + \Mper_{53421}.
\end{split}
\end{equation}

\end{example}

\subsection{Duality}

Recall that the coproduct is given by

\begin{equation}
\label{coprodM}
\Delta \Mper_\sigma := \sum_{(\alpha,\beta) ; \alpha\sconc \beta=\sigma} 
\Mper_{\alpha} \tensor \Mper_{\beta}.
\end{equation}
As in Section~\ref{eqsym}, this implies
\begin{proposition}
If $(\Sper^\sigma)$ denotes the dual basis of $(\Mper_\sigma)$, the graded
dual $\SGSym := \SGQSym^*$ is free over the set
\begin{equation}
\{ \Sper^\alpha \,|\, \alpha \text{\ connected} \}.
\end{equation}
\end{proposition}

Indeed, Equation~(\ref{coprodM}) is equivalent to
\begin{equation}
\Sper^\alpha \Sper^\beta = \Sper^{\alpha\sconc\beta}.
\end{equation}
\qed

Note that $\SGSym$ is both a subalgebra and a quotient of $\ESym$, since
$\SGQSym$ is both a quotient and a subalgebra of $\EQSym$.

Now, as before, $\SGSym$ being a graded connected cocommutative Hopf algebra,
it follows from the Milnor-Moore theorem that
\begin{equation}
\SGSym = U(L),
\end{equation}
where $L$ is the Lie algebra of its primitive elements.

The same argument as in Section~\ref{eqsym} proves
\begin{theorem}
The graded Lie algebra $L$ of primitive elements of $\SGSym$ is free over a
set indexed by connected permutations.
\qed
\end{theorem}

\begin{corollary}
$\SGSym$ is isomorphic to $H_O$, the Grossman-Larson Hopf algebra of
heap-ordered trees \cite{GL}.
\qed
\end{corollary}

According to~\cite{AS}, $\SGQSym$ ($=\SGSym^*$) is therefore isomorphic to the
quotient of $\FQSym$ by its coradical filtration.

%
\subsection{Cyclic tensors and $\SGQSym$}

For a vector space $V$, let $\Gamma^n V$ be the subspace of $V^{\otimes n}$
spanned by \emph{cyclic tensors}, \emph{i.e.}, sums of the form
\begin{equation}
\sum_{k=0}^{n-1} (v_1\otimes \cdots\otimes v_k) \gamma^k,
\end{equation}
where $\gamma$ is the cycle $(1,2,\ldots,n)$, the right action of permutations
on tensors being as usual
\begin{equation}
(v_1\otimes \cdots\otimes v_k) \sigma =
v_{\sigma(1)}\otimes \cdots \otimes v_{\sigma(n)}.
\end{equation}

Clearly, $\Gamma^n V$ is stable under the action of $GL(V)$, and its character
is the symmetric function ``cyclic character''~\cite{ST,LST}:
\begin{equation}
l_n^{(0)} = \frac{1}{n} \sum_{d|n} \phi(d) p_d^{n/d},
\end{equation}
where $\phi$ is Euler's function.

Let $L_n^{(0)}$ be the subspace of $\CC \SG_n$ spanned by $n$-cycles. This is
a submodule of $\CC\SG_n$ for the conjugation action
$\rho_\tau(\sigma)=\tau\sigma\tau^{-1}$ with Frobenius characteristic
$l_n^{(0)}$.
Then one can define the analytic functor $\Gamma$~\cite{Joy,Mcd}:
\begin{equation}
\Gamma (V) = \sum_{n\geq0} V^{\otimes n} \otimes_{\CC\SG_n} L_n^{(0)}.
\end{equation}

Let $\overline{\Gamma} (V) = \bigoplus_{n\geq1} \Gamma^n(V)$. Its symmetric
algebra $H(V)=S(\overline{\Gamma} (V))$ can be endowed with a Hopf algebra
structure, by declaring the elements of $\overline{\Gamma} (V)$ primitive.

As an analytic functor, $V\mapsto S(\overline{\Gamma} (V))$ corresponds to the
sequence of $\SG_n$-modules $M_n = \CC \SG_n$ endowed with the conjugation
action, that is,
\begin{equation}
S(\overline{\Gamma} (V)) = \bigoplus_{n\geq0} V^{\otimes n}\otimes_{\CC\SG_n}
                          M_n,
\end{equation}
so that elements of $H_n(V)$ can be identified with symbols
$\qbin{w}{\sigma}$ with $w\in V^{\otimes n}$ and $\sigma\in\SG_n$ subject to
the equivalences
\begin{equation}
\qbin{w \tau^{-1}}{\tau\sigma\tau^{-1}} \equiv \qbin{w}{\sigma}.
\end{equation}
Let $A=\{a_n\,|\, n\geq 1\}$ be an infinite linearly ordered alphabet, and
$V=\CC A$. We identify a tensor product of letters
$a_{i_1}\otimes\cdots\otimes a_{i_n}$
with the corresponding word $w=a_{i_1}\ldots a_{i_n}$ and denote by $(w)$ the
circular class of $w$. A basis of $H_n$ is then formed by the commutative
products
\begin{equation}
\underline{m} = (w_1)\cdots(w_p)
\end{equation}
of circular words, with $|w_1|+\cdots+|w_p|=n$.

With such a basis element, we can associate a permutation by the following
standardization process. Fix a total order on circular words, for example the
lexicographic order on minimal representatives. Write $\underline{m}$ as a
non-decreasing product
\begin{equation}
\underline{m} = (w_1)\cdots(w_p) \text{ with } (w_1)\leq
(w_2)\leq\cdots\leq(w_p),
\end{equation}
and compute the ordinary standardization $\sigma'$ of the word
$w=w_1\cdots w_p$. Then, $\sigma$ is the permutation obtained by parenthesing
the word $\sigma'$ like $\underline{m}$ and interpreting the factors as
cycles. For example, if
\begin{equation}
\begin{split}
\underline{m} &= (cba)(aba)(ac)(ba) = (aab)(ab)(ac)(acb) \\
w &= aababacacb \\
\sigma' &= 12637495\,10\,8 \\
\sigma &= (126)(37)(49)(5\,10\,8) \\
\sigma &= (2,6,7,9,10,1,3,5,4,8) \\
\end{split}
\end{equation}
We set $\sigma = \cstd(\underline{m})$ and define it as the \emph{circular
standardized} of $m$.

Let $H_\sigma(V)$ be the subspace of $H_n(V)$ spanned by those
$\underline{m}$ such that $\cstd(\underline{m})=\sigma$, and let
$\pi_\sigma : H(V)\to H_\sigma(V)$ be the projector associated with the direct
sum decomposition
\begin{equation}
H(V) = \bigoplus_{n\geq0} \bigoplus_{\sigma\in\SG_n} H_\sigma(V).
\end{equation}

Computing the convolution of such projectors then yields the following
\begin{theorem}
The $\pi_\sigma$ span a subalgebra of the convolution algebra
$\End^{\gr} H(V)$, isomorphic to $\SGQSym$ via $\pi_\sigma\mapsto M_\sigma$.
\qed
\end{theorem}

\subsection{Subalgebras of $\SGQSym$}

For a permutation $\sigma\in\SG_n$, let $\csupp(\sigma)$ be the partition
$\pi$ of the set $[n]$ whose blocks are the supports of the cycles of
$\sigma$. The sums
\begin{equation}
\upi_{\pi} := \sum_{\csupp(\sigma)=\pi} \Mper_\sigma
\end{equation}
span a Hopf subalgebra $\PiQSym$ of $\SGQSym$, which, as we shall see in the
next section, is isomorphic to the dual of the Hopf algebra of symmetric
functions in noncommuting variables (such as in~\cite{SR,BRRZ}, not to be
confused with $\NCSF$).

\medskip
One can also embed $\QSym$ into $\PiQSym$: take as total ordering on finite
sets of integers $\{i_1<\cdots<i_r\}$ the lexicographic order on the words
$(i_1,\ldots,i_r)$. Then, any set partition $\pi$ of $[n]$ has a canonical
representative $B$ as a non-decreasing sequence of blocks $(B_1\leq B_2\leq
\cdots\leq B_r)$. Let $I=b(\pi)$ be the composition $(|B_1|,\ldots,|B_r|)$ of
$n$.
The sums
\begin{equation}
\uq_{I} := \sum_{b(\pi)=I} \upi_\pi
 = \sum_{b(\sigma)=\lambda} \Mper_\sigma
\end{equation}
where $b(\sigma)$ denotes the ordered cycle type of $\sigma$, span a Hopf
subalgebra of $\PiQSym$ and $\SGQSym$, which is isomorphic to $\QSym$.

\medskip
Furthemore, if we denote by $\Lambda(I)$ the partition associated with a
composition $I$ by sorting $I$ and by $\Lambda(\pi)$ the partition $\lambda$
whose parts are the sizes of the blocks of $\pi$, the sums
\begin{equation}
\ul_{\lambda} := \sum_{\Lambda(I)=\lambda} \uq_{I}
 = \sum_{\Lambda(\pi)=\lambda} \upi_\pi
 = \sum_{C(\sigma)=\lambda} \Mper_\sigma
\end{equation}
where $C(\sigma)$ denotes the cycle type of $\sigma$, span a Hopf subalgebra
of $\QSym$, $\PiQSym$, and $\SGQSym$, which is isomorphic to $Sym$ (ordinary
symmetric functions).

An explicit Hopf embedding of $Sym$ into $\SGQSym$ is given by
\begin{equation}
j: p_\lambda^* \to \ul_\lambda
\end{equation}
where $p_\lambda^* = \frac{p_\lambda}{z_\lambda}$ is the adjoint basis of
products of power sums.
The images of the usual generators of $Sym$ under this embedding have simple
expressions in terms of the infinite matrix $X=(x_{ij})_{i,j\leq1}$:
\begin{equation}
j(p_n) = \tr(X^n)
\end{equation}
which implies that $j(e_n)$ is the sum of the diagonal minors of order $n$ of
$X$:
\begin{equation}
j(e_n) = \sum_{i_1<\cdots<i_n} \sum_{\sigma\in\SG_n}
         \varepsilon(\sigma) x_{i_1 i_{\sigma(1)}}\ldots x_{i_n i_{\sigma(n)}}
\end{equation}
and $j(h_n)$ is the sum of the same minors of the permanent
\begin{equation}
j(h_n) = \sum_{i_1<\cdots<i_n} \sum_{\sigma\in\SG_n}
         x_{i_1 i_{\sigma(1)}} \ldots x_{i_n i_{\sigma(n)}}.
\end{equation}
More generally, the sum of the diagonal immanants of type $\lambda$ gives
\begin{equation}
j(s_\lambda) =
\sum_{i_1<\cdots<i_n} \sum_{\sigma\in\SG_n}
\chi^\lambda(\sigma) x_{i_1 i_{\sigma(1)}} \ldots x_{i_n i_{\sigma(n)}}.
\end{equation}

\medskip
Finally, one can check that the $M_\sigma$ with $\sigma$ involutive span a
Hopf subalgebra of $\SGQSym$. Since the number of involutions of $\SG_n$ is
equal to the number of standard Young tableaux of size $n$, this algebra can
be regarded as a commutative version of $\FSym$.
Notice that this version is also isomorphic to the image of $\SGQSym$ in the
quotient of $R/{\mathcal J}$ by the relations
\begin{equation}
x_{ij}x_{jk}=0, \text{\ for all $i\not=k$.}
\end{equation}
This construction generalizes to the algebras built on permutations of
arbitrary given order.

\section{Structure of $\SGSym$}

\subsection{A realization of $\SGSym$}

In the previous section, we have built a commutative algebra of permutations
from explicit polynomials on a set of auxiliary variables $x_{i\,j}$.
One may ask whether its non-commutative dual admits a similar realization in
terms of non-commuting variables $a_{i\,j}$.

We shall find such a realization, in a somewhat indirect way, by first
building from scratch a Hopf algebra of permutations
$\PhiSym \subset \KK\,\langle\, a_{i\,j} \,|\, i,j\geq 1\,\rangle$,
whose operations can be described in terms of the cycle structure of
permutations.
Its coproduct turns out to be cocommutative, and the isomorphism with
$\SGSym$ follows as above from the Milnor-Moore theorem.

Let $\{a_{i\,j}, i,j\geq1\}$ be an infinite set of non-commuting
indeterminates. We use the biword notation
\begin{equation}
a_{i\,j} \equiv \ncbinomial{i}{j}, \quad
\ncbinomial{i_1}{j_1} \cdots \ncbinomial{i_n}{j_n} \equiv
\ncbinomial{i_1\ldots i_n}{j_1\ldots j_n}
\end{equation}

Let $\sigma\in\SG_n$ and let $(c_1,\ldots,c_k)$ be a decomposition of $\sigma$
into disjoint cycles. With any cycle, one associates its \emph{cycle words},
that is, the words obtained by reading the successive images of any
element of the cycle.
For example, the cycle words associated with the cycle $(3142)$ are $1342,
2134, 3421, 4213$.

We now define
\begin{equation}
\myphi_\sigma := \sum_{x,a} \ncbinomial{x}{a},
\end{equation}
where the sum runs over all words $x$ such that $x_i=x_j$ iff $i$ and
$j$ belong to the same cycle of $\sigma$, and such that the standardized word
of the subword of $a$ whose letter positions belong to cycle $c_l$ is equal to
the inverse of the standardized word of a cycle word of $c_l$.

\begin{example}
\begin{equation}
\myphi_{12} = \sum_{x\not= y} \ncbinomial{x\ y}{a\ b}.
\end{equation}
\begin{equation}
\myphi_{41352} =
  \sum_{x\not=y ; \Std(abde)^{-1}=1342, 3421, 4213, \text{\rm or\ } 2134}
  {\ncbinomial{x\ x\ y\ x\ x}{a\ b\ c\ d\ e}}.
\end{equation}
\end{example}

\begin{theorem}
The $\myphi_\sigma$ span a subalgebra $\PhiSym$ of
$\KK\,\langle\, a_{ij} \,|\, i,j\geq1 \,\rangle$.
More precisely, there exist non-negative integers
$g_{\alpha,\beta}^{\sigma}$ ($0$ or $1$) such that
$\myphi_\alpha \myphi_\beta = \sum g_{\alpha,\beta}^{\sigma} \myphi_\sigma$.
\medskip

$\PhiSym$ is free over the set
\begin{equation}
\{\myphi_\alpha \,|\, \alpha \text{\ connected} \}.
\end{equation}
\end{theorem}

To give the precise expression of the product $\myphi_\alpha\myphi_\beta$, we 
first need to define two operations on cycles.

The first operation is just the circular shuffle on disjoint cycles: if $c'_1$
and $c''_1$ are two disjoint cycles, their \emph{cyclic shuffle}
$c'_1 \picyc c''_1$ is the set of cycles $c_1$ such that their cycle words are
obtained by applying the usual shuffle on the cycle words of $c'_1$ and
$c''_1$. This definition makes sense because a shuffle of cycle words
associated with two words on disjoint alphabets splits as a union of cyclic
classes.

For example, the cyclic shuffle $(132)\picyc(45)$ gives the set of cycles
\begin{equation}
\begin{split}
\{ (13245), (13425), (13452), (14325), (14352), (14532), \\
   (13254), (13524), (13542), (15324), (15342), (15432) \}.
\end{split}
\end{equation}
These cycles correspond to the following list of permutations which are
those appearing in Equation~(\ref{312-54b}), except for the first one which
will be found later:
\begin{equation}
\begin{split}
\{ &\ 34251,\ 35421,\ 31452,\ 45231,\ 41532,\ 41253,  \\
   &\ 35214,\ 34512,\ 31524,\ 54213,\ 51423,\ 51234 \}.
\end{split}
\end{equation}

Let us now define an operation on two sets $C_1$ and $C_2$ of disjoint cycles.
We call \emph{matching} a list of all those cycles, some of the cycles being
paired, always one of $C_1$ with one of $C_2$. The cycles remaining alone are
considered to be associated with the empty cycle. With all matchings associate
the set of sets of cycles obtained by the product $\picyc$ of any pair of
cycles.
The union of those sets of cycles is denoted by $C_1\picycgen C_2$.

For example, the matchings corresponding to
$C_1=\{(1),(2)\}$ and $C_2=\{(3),(4)\}$ are:
\begin{equation}
\label{c1c2a}
\begin{split}
&\{(1)\}\{(2)\}\{(3)\}\{(4)\},\
\{(1)\}\{(2),(3)\}\{(4)\},\ 
\{(1)\}\{(2),(4)\}\{(3)\}, \\
&\{(1),(3)\}\{(2)\}\{(4)\},\
\{(1),(3)\}\{(2),(4)\},\
\{(1),(4)\}\{(2)\}\{(3)\},\
\{(1),(4)\}\{(2),(3)\},
\end{split}
\end{equation}
and the product $C_1\picycgen C_2$ is then
\begin{equation}
\label{c1c2b}
\begin{split}
&\{(1),(2),(3),(4)\},\ \
\{(1),(23),(4)\},\ \
\{(1),(24),(3)\}, \\
&\{(13),(2),(4)\},\ \
\{(13),(24)\},\ \
\{(14),(2),(3)\},\ \
\{(14),(23)\}.
\end{split}
\end{equation}

Remark that this calculation is identical with the Wick formula in quantum
field theory (see~\cite{BO} for an explanation of this coincidence).

We are now in a position to describe the product $\myphi_\sigma \myphi_\tau$:
let $C_1$ be the cycle decomposition of $\sigma$ and $C_2$ be the cycle
decomposition of $\tau$, shifted by the size of $\sigma$. Then the
permutations indexing the elements appearing in the product
$\myphi_\sigma \myphi_\tau$ are the permutations whose cycle decompositions
belong to $C_1\picycgen C_2$.

For example, with $\sigma=\tau=12$, one finds that $C_1=\{(1),(2)\}$ and
$C_2=\{(3),(4)\}$. It is then easy to check that one goes from
Equation~(\ref{c1c2b}) to Equation~(\ref{12-12b}) by computing the
corresponding permutations.

\begin{example}
{\rm
\begin{equation}
\label{12-43b}
\myphi_{12} \myphi_{21} = \myphi_{1243} + \myphi_{1342} + \myphi_{1423} +
\myphi_{3241} + \myphi_{4213}.
\end{equation}
\begin{equation}
\label{12-12b}
\myphi_{12} \myphi_{12} = \myphi_{1234} + \myphi_{1324} + \myphi_{1432} +
\myphi_{3214} + \myphi_{3412} + \myphi_{4231} + \myphi_{4321}.
\end{equation}
\begin{equation}
\myphi_{1} \myphi_{4312} = \myphi_{15423} + \myphi_{25413} + \myphi_{35421} +
\myphi_{45123} + \myphi_{51423}.
\end{equation}
\begin{equation}
\begin{split}
\label{312-54b}
\myphi_{312} \myphi_{21} &=\ 
\myphi_{31254} + \myphi_{31452} + \myphi_{31524} + \myphi_{34251} +
\myphi_{34512} + \myphi_{35214} + \myphi_{35421} \\
& +\ \myphi_{41253} + \myphi_{41532} +  \myphi_{45231} + \myphi_{51234} +
\myphi_{51423} + \myphi_{54213}.
\end{split}
\end{equation}
}
\end{example}

Let us recall a rather general recipe to obtain the coproduct of a
combinatorial Hopf algebra from a realization in terms of words on an ordered
alphabet $X$.
Assume that $X$ is the ordered sum of two mutually commuting alphabets
$X'$ and $X''$.  Then define the coproduct as $\Delta(F)=F(X'\dot{+}X'')$,
identifying $F'\tensor F''$ with $F'(X')F''(X'')$~\cite{NCSF6, NTpfb}.

There are many different ways to define a coproduct on $\PhiSym$ compatible
with the realization since there are many ways to order an alphabet of
biletters: order the letters of the first alphabet, order the letters of the
second alphabet, or order lexicographically with respect to one alphabet and
then to the second.

In the sequel, we only consider the coproduct obtained by ordering the
biletters with respect to the first alphabet. Thanks to the definition of the
$\myphi$, it is easy to see that it corresponds to the unshuffling of the
cycles of a permutation:

\begin{equation}
\Delta\myphi_\sigma := \sum_{(\alpha,\beta)}
\myphi_\alpha \tensor \myphi_\beta,
\end{equation}
where the sum is taken over all pairs of permutations $(\alpha,\beta)$ such
that $\alpha$ is obtained by standardizing any subset of cycles of $\sigma$,
and $\beta$ by standardizing the complementary subset of cycles.

\begin{example}
\begin{equation}
\Delta\myphi_{12} = \myphi_{12}\tensor 1 + 2 \myphi_{1}\tensor\myphi_{1}
+1 \tensor \myphi_{12}.
\end{equation}
\begin{equation}
\Delta\myphi_{312} = \myphi_{312}\tensor 1 + 1 \tensor \myphi_{312}.
\end{equation}
\begin{equation}
\Delta\myphi_{4231} = \myphi_{4231}\tensor 1 + 2 \myphi_{321}\tensor\myphi_{1}
+ \myphi_{21}\tensor \myphi_{12} + \myphi_{12}\tensor \myphi_{21}
+ 2 \myphi_{1}\tensor \myphi_{321} + 1\tensor \myphi_{4231}.
\end{equation}
\end{example}

The next theorem can be easily proved on the realization.

\begin{theorem}
$\Delta$ is an algebra morphism, so that $\PhiSym$ is a graded
bialgebra (for the grading $\deg\myphi_\sigma = n$ if $\sigma\in\SG_n$).
Moreover, $\Delta$ is cocommutative.
\qed
\end{theorem}

The same reasoning as in Section~\ref{sgqsym} shows that

\begin{theorem}
$\SGSym$ and $\PhiSym$ are isomorphic as Hopf algebras.
\qed
\end{theorem}

To get the explicit change of basis going from $\myphi$ to $\Sper$, let us
first recall that a \emph{connected permutation} is a permutation $\sigma$
such that $\sigma([1,k])\not = [1,k]$ for any $k\in [1,n-1]$. Any permutation
$\sigma$ has a unique maximal factorization
$\sigma=\sigma_1\bullet\cdots\bullet\sigma_r$ into connected permutations.
We then define
\begin{equation}
\label{sprime}
\Sper'_\sigma := \myphi_{\sigma_1}\cdots\myphi_{\sigma_r}.
\end{equation}

First remark that the $\Sper'$ form a basis of $\PhiSym$.
It is easy to check that the $\Sper'$ is a multiplicative basis with product
given by shifted concatenation of permutations, so that they multiply as the
$\Sper$ do. Moreover, the coproduct of $\Sper'_\sigma$ is the same as for
$\myphi_\sigma$, so the same as for $\Sper^\sigma$.
So both bases $\Sper$ and $\Sper'$ have same product and same coproduct. 

This proves that the linear map $ \Sper^\sigma \mapsto \Sper'_\sigma $
realizes the Hopf isomorphism between $\SGSym$ and $\PhiSym$.
There is another natural isomorphism: define
\begin{equation}
\label{ssec}
\Sper''_\sigma := \sum_{x,a} \ncbinomial{x}{a},
\end{equation}
where the sum runs over all words $x$ such that $x_i=x_j$ if (but \emph{not}
only if) $i$ and $j$ belong the same cycle of $\sigma$ and such that the
standardized word of the subword of $a$ consisting of the indices of cycle
$c_l$ is equal to the inverse of the standardized word of a cycle word of
$c_l$.

The fact that both bases $\Sper'$ and $\Sper''$ have same product and
coproduct simply comes from the fact that if $(c_1)\cdots(c_p)$ is the cycle
decomposition of $\sigma$,

\begin{equation}
\Sper''_\sigma = \sum_{(c);(c)\in
                 (c_1)\picycgen(c_2)\picycgen\cdots\picycgen(c_p)}
                     \myphi_{(c)}.
\end{equation}

For example,
\begin{equation}
\begin{split}
\Sper''_{2431} = \Sper''_{(124)(3)}  & =
\myphi_{(124)(3)} + \myphi_{(1423)} + \myphi_{(1234)} + \myphi_{(1324)}\\
& =\myphi_{2431}  + \myphi_{4312}   + \myphi_{2341}   + \myphi_{3421}.
\end{split}
\end{equation}

\subsection{Quotients of $\PhiSym$}

Let $I$ be the ideal of $\PhiSym$ generated by the differences
\begin{equation}
\phi_{\sigma} - \phi_{\tau}
\end{equation}
where $\sigma$ and $\tau$ have the same cycle type.

The definitions of its product and coproduct directly imply that $I$ is a Hopf
ideal.
Since the cycle types are parametrized by integer partitions, the quotient
$\PhiSym/I$ has a basis $Y_\lambda$, corresponding to the class of
$\phi_{\sigma}$, where $\sigma$ has $\lambda$ as cycle type.

From Equations~(\ref{12-43b})-(\ref{312-54b}), one finds:

\begin{example}
\begin{equation}
Y_{11} Y_{2} = Y_{211} + 4 Y_{31}, \qquad
Y_{11}^2 = Y_{1111} + 2 Y_{22} + 4 Y_{211}.
\end{equation}
\begin{equation}
Y_1 Y_4 = Y_{41} + 4 Y_5, \qquad Y_{3} Y_{2} = Y_{32} + 12 Y_5.
\end{equation}
\end{example}

\begin{theorem}
$\PhiSym/I$ is isomorphic to $Sym$, the Hopf algebra of ordinary symmetric
functions,
\medskip

If one writes
$\lambda=(\lambda_1,\ldots,\lambda_p)=(1^{m_1},\ldots,k^{m_k})$,
an explicit isomorphism is given by
\begin{equation}
Y_\lambda \mapsto 
\frac{\prod_{i=1}^{k} m_i!}{\prod_{j=1}^p (\lambda_j-1)!} m_\lambda.
\end{equation}
\qed
\end{theorem}

\section{Parking functions and trees}

\subsection{A commutative algebra of parking functions}

It is also possible to build a commutative pendant of the Hopf algebra of
parking functions introduced in~\cite{NT1}:
let $\PF_n$ be the set of parking functions of length $n$.
For $\park\in\PF_n$, set, as before

\begin{equation}
\Mpa_\park := \sum_{i_1<\cdots < i_n} x_{i_1\ i_{\park(1)}} \cdots
x_{i_n\ i_{\park(n)}}.
\end{equation}

Then, once more, the $\Mpa_\park$ form a linear basis of a $\ZZ$-subalgebra
$\CPQSym$ of $\EQSym$, which is also a sub-coalgebra if one defines the
coproduct in the usual way, that is, from special cuts in graphs
(see~\cite{NTT} for more details).

\begin{example}
\begin{equation}
\Mpa_{1} \Mpa_{11} = \Mpa_{122} + \Mpa_{121} + \Mpa_{113}.
\end{equation}
\begin{equation}
\Mpa_{1} \Mpa_{221} = \Mpa_{1332} + \Mpa_{3231} + \Mpa_{2231} + \Mpa_{2214}.
\end{equation}
\begin{equation}
\Mpa_{12} \Mpa_{21} = \Mpa_{1243} + \Mpa_{1432} + \Mpa_{4231} + \Mpa_{1324} +
\Mpa_{3214} + \Mpa_{2134}.
\end{equation}
\begin{equation}
\Delta \Mpa_{525124}  = \Mpa_{525124} \tensor 1 + 1\tensor \Mpa_{525124}.
\end{equation}
\begin{equation}
\Delta \Mpa_{4131166} = \Mpa_{4131166} \tensor 1 +
                        \Mpa_{41311}\tensor \Mpa_{11} +
                        1 \tensor \Mpa_{4131166}.
\end{equation}
\end{example}

The main interest of the non-commutative and non-cocommutative Hopf algebra of
parking functions defined in~\cite{NT1} was that it naturally led to two
algebras of trees. We obtained a cocommutative Hopf algebra of planar binary
trees by summing over the distinct permutations of parking functions, and an
algebra of planar trees by summing over hypoplactic classes.

We shall now investigate whether similar constructions can be found for the
commutative version $\CPQSym$.

\subsection{From labelled to unlabelled parking graphs}

A first construction, which can always be done for Hopf algebras of labelled
graphs is to build a subalgebra by summing over labellings.
Notice that this subalgebra is the same as the subalgebra obtained by summing
endofunctions graphs over their labellings.

The dimension of this subalgebra in degree $n$ is equal to the number of
unlabelled parking graphs $1$, $1$, $3$, $7$, $19$, $47$, $\ldots$
For example, here are the $7$ unlabelled parking graphs of size $3$ (to be
compared with the $16$ parking functions):

\begin{figure}[ht]
\epsfig{file=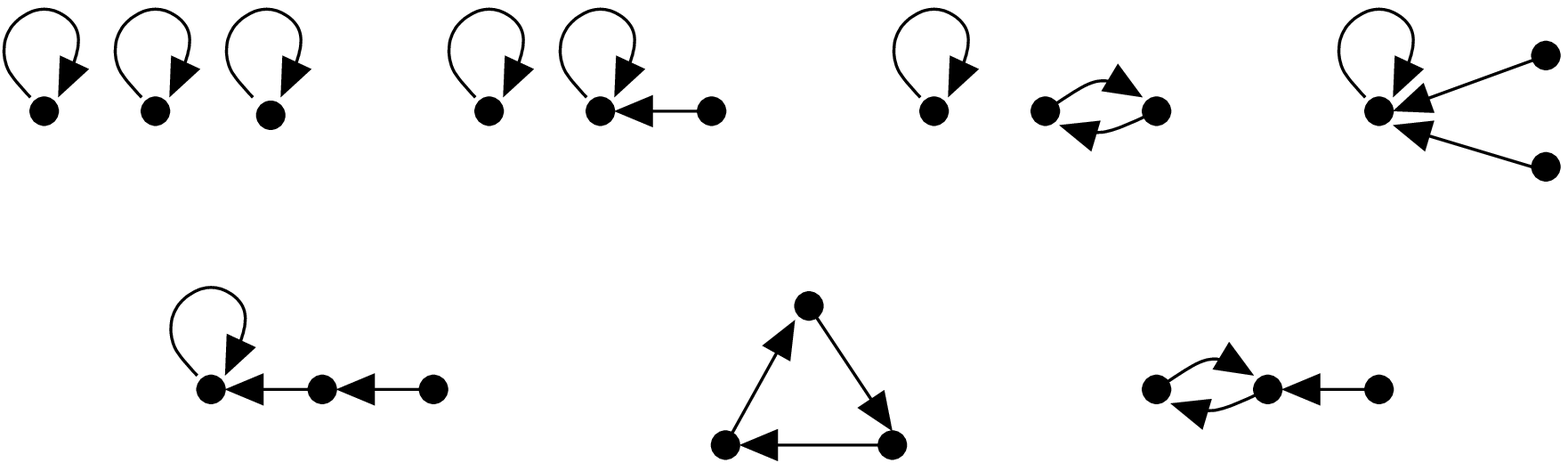,width=12cm}
\end{figure}

The product of two such unlabelled graphs is the concatenation of graphs and
the coproduct of an unlabelled graph is the unshuffle of its connected
subgraphs. So this algebra is isomorphic to the polynomial algebra on
generators indexed by connected parking graphs.

\subsection{Binary trees and nondecreasing parking functions}

One can easily check that in $\CPQSym$, summing over parking functions having
the same reordering does not lead to a subalgebra. However, if we denote by
$I$ the subspace of $\CPQSym$ spanned by the $\Mpa_\park$ where $\park$ is not
equal to its nondecreasing reordering, it turns out that $I$ is an ideal and a
coideal, and $\CCQSym := \CPQSym/I$ is therefore a commutative Hopf algebra
with basis given by the classes $\Mpa_{\parkc}:=\Mpa_{\overline\park}$
labelled by nondecreasing parking functions.

Notice that $\CCQSym$ is also isomorphic to the image of $\CPQSym$ in the
quotient of $R/{\mathcal J}$ by the relations
\begin{equation}
x_{ij}x_{kl}=0, \text{\ for all $i<k$ and $j>l$.}
\end{equation}

The dual basis of $\Mpa_{\parkc}$ is
\begin{equation}
S^{\parkc} := \sum_{\park} S^\park,
\end{equation}
where the sum is taken over all permutations of $\parkc$.

The dual $\CCQSym^*$ is free over the set $S^\parkc$, where $\parkc$ runs over
connected nondecreasing parking functions. So if one denotes by $\CQSym$ the
Catalan algebra defined in~\cite{NT1}, 
the usual Milnor-Moore argument then shows that
\begin{equation}
\CQSym \sim \CCQSym^* ,\qquad\qquad \CCQSym \sim \CQSym^*.
\end{equation}

\subsection{From nondecreasing parking functions to rooted forests}

Nondecreasing parking functions correspond to parking graphs of a particular
type: namely, rooted forests with a particular labelling (it corresponds to
nondecreasing maps), the root being given by the loops in each connected
component.

Taking sums over the allowed labellings of a given rooted forest, we get that
the
\begin{equation}
M_F := \sum_{supp(\parkc)=F} M_\parkc,
\end{equation}
span a commutative Hopf algebra of rooted forests, which is likely to coincide
with the quotient of the Connes-Kreimer algebra~\cite{CK} by its coradical
filtration~\cite{AS}.

\section{Quantum versions}

\subsection{Quantum quasi-symmetric functions}

When several Hopf algebra structures can be defined on the same class of
combinatorial objects, it is tempting to try to interpolate between them.

This can be done, for example with compositions: the algebra of quantum
quasi-symmetric functions $\QSym_q$~\cite{TU,NCSF6} interpolates between
quasi-symmetric functions and non-commutative symmetric functions.

However, the natural structure on $\QSym_q$ is not exactly that of a Hopf
algebra but rather of a \emph{twisted Hopf algebra}~\cite{LZ}.

Recall that the coproduct of $\QSym(X)$ amounts to replace $X$ by the ordered
sum $X' \dot{+} X''$ of two isomorphic and mutually commuting alphabets. On
the other hand, $\QSym_q$ can be realized by means of an alphabet of
$q$-commuting letters
\begin{equation}
x_jx_i = q x_i x_j, \text{ for } j>i.
\end{equation}
Hence, if we define a coproduct on $\QSym_q$ by
\begin{equation}
\Delta_q f(X) = f(X'\dot+X''),
\end{equation}
with $X'$ and $X''$ $q$-commuting with each other, it will be an algebra
morphism
\begin{equation}
\QSym_q \to \QSym_q(X'\dot+X'') = \QSym_q \otimes_{\chi} \QSym_q
\end{equation}
for the \emph{twisted tensor product}
\begin{equation}
(a\otimes b) \cdot (a'\otimes b') = \chi(b,a') (aa' \otimes bb'),
\end{equation}
where
\begin{equation}
\label{chi}
\chi(b,a') = q^{\deg(b)\cdot \deg(a')}
\end{equation}
for homogeneous elements $b$ and $a'$.

It is easily checked that $\Delta_q$ is actually given by the same formula as
the usual coproduct of $\QSym$, that is
\begin{equation}
\Delta_q M_I = \sum_{H\cdot K=I} M_H\otimes M_K.
\end{equation}

The dual twisted Hopf algebra, denoted by $\NCSF_q$, is isomorphic to $\NCSF$
as a algebra. If we denote by $S^I$ the dual basis of $M_I$, $S^I
S^J=S^{I\cdot J}$, and $\NCSF_q$ is freely generated by the $S^{(n)}=S_n$,
whose coproduct is
\begin{equation}
\Delta_q S_n = \sum_{i+j=n} q^{ij} S_i\otimes S_j.
\end{equation}

As above, $\Delta_q$ is an algebra morphism
\begin{equation}
\NCSF_q \to \NCSF_q \otimes_{\chi} \NCSF_q,
\end{equation}
where $\chi$ is again defined by Equation~(\ref{chi}).

\subsection{Quantum free quasi-symmetric functions}

The previous constructions can be lifted to $\FQSym$.
Recall that $\phi(\F_\sigma) = q^{l(\sigma)} F_{c(\sigma)}$ is an algebra
homomorphism $\FQSym \to \QSym_q$, which is in fact induced by the
specialization $\phi(a_i)=x_i$ of the underlying free variables $a_i$ to
$q$-commuting variables $x_i$.

The coproduct of $\FQSym$ is also defined by
\begin{equation}
\Delta \F(A) = \F(A' \dot+ A''),
\end{equation}
where $A'$ and $A''$ are two mutually commuting copies of $A$~\cite{NCSF6}.
If instead one sets $a''a' = qa' a''$, one obtains again a twisted Hopf
algebra structure $\FQSym_q$ on $\FQSym$, for which $\phi$ is a
homomorphism.

\medskip
\begin{theorem}
~
Let $A'$ and $A''$ be $q$-commuting copies of the ordered alphabet $A$,
\emph{i.e.}, $a''a'=qa'a''$ for $a'\in A'$ and $a''\in A''$.
Then, the coproduct
\begin{equation}
\Delta_q f = f(A'\dot+ A'')
\end{equation}
defines a twisted Hopf algebra structure.
It is explictly given in the basis $\F_\sigma$ by
\begin{equation}
\label{deltaq}
\Delta_q \F_\sigma = \sum_{\alpha\sconc\beta=\sigma}
                     q^{\inv(\alpha,\beta)} \F_\alpha\otimes \F_\beta
\end{equation}
where $\inv(\alpha,\beta)$ is the number of inversions of $\sigma$ with one
element in $\alpha$ and the other in $\beta$.
\medskip

More precisely, $\Delta_q$ is an algebra morphism with values in the
twisted tensor product of graded algebras $\FQSym\otimes_\chi \FQSym$
where $(a\otimes_\chi b)(a'\otimes_\chi b') = \chi(b,a')
(aa'\otimes_\chi bb')$ and $\chi(b,a')=q^{\deg(b).\deg(a')}$ for homogeneous
elements $b$, $a'$.
\medskip

The map $\phi: \FQSym_q \to \QSym_q$ defined by
\begin{equation}
\phi(\F_\sigma) = q^{l(\sigma)}  F_{c(\sigma)}
\end{equation}
is a morphism of twisted Hopf algebras.
\qed
\end{theorem}

\begin{example}
\begin{equation}
\Delta_q \F_{2431}  = \F_{2431}\otimes 1 + q^3 \F_{132}\otimes \F_{1}
+ q^3 \F_{12}\otimes \F_{21} + q\F_1\otimes\F_{321} + 1\otimes \F_{2431}.
\end{equation}
\begin{equation}
\Delta_q \F_{3421} =  \F_{3421}\otimes 1 + q^3 \F_{231}\otimes \F_{1}
+ q^4 \F_{12}\otimes \F_{21} + q^2\F_1\otimes\F_{321} + 1\otimes \F_{3421}.
\end{equation}
\begin{equation}
\Delta_q \F_{21} = \F_{21}\otimes 1 + q \F_{1}\otimes\F_{1} + 1\otimes\F_{21}.
\end{equation}
\begin{equation}
\begin{split}
(\Delta_q \F_{21})(\Delta_q\F_{1}) =
(\F_{213}+\F_{231}+\F_{321}) \otimes 1
+ ( \F_{21} + q^2(\F_{12}+\F_{21}))\otimes \F_{1} \\
+ \F_{1} \otimes ( q^2\F_{21} + q(\F_{12}+\F_{21}))
+1\otimes (\F_{213}+\F_{231}+\F_{321}).
\end{split}
\end{equation}
\begin{equation}
\Delta_q \F_{213} = \F_{213}\otimes 1 + \F_{21}\otimes \F_{1}
+ q \F_{1}\otimes \F_{12} + 1\otimes \F_{213}.
\end{equation}
\begin{equation}
\Delta_q \F_{231} = \F_{231}\otimes 1 + q^2 \F_{12}\otimes \F_{1}
+ q \F_{1}\otimes \F_{21} + 1\otimes \F_{231}.
\end{equation}
\begin{equation}
\Delta_q \F_{321} = \F_{321}\otimes 1 + q^2 \F_{21}\otimes \F_{1}
+ q^2 \F_{1}\otimes \F_{21} + 1\otimes \F_{321}.
\end{equation}
\end{example}

Finally, one can also define a one-parameter family of ordinary Hopf algebra
structures on $\FQSym$, by restricting formula~(\ref{deltaq}) for $\Delta_q$
to connected permutations $\sigma$, and requiring that $\Delta_q$ be an
algebra homomorphism. Then, for $q=0$, $\Delta_q$ becomes cocommutative, and
it is easily shown that the resulting Hopf algebra is isomorphic to $\SGSym$.

However, it follows from~\cite{NCSF6} that for generic $q$, the Hopf algebras
defined in this way are all isomorphic to $\FQSym$.
This suggests to interpret $\FQSym$ as a kind of quantum group: it would be
the generic element of a quantum deformation of the enveloping algebra
$\SGSym=U(L)$.
Similar considerations apply to various examples, in particular to the
Loday-Ronco algebra $\PBT$, whose commutative version obtained in~\cite{NT1}
can be quantized in the same way as $\QSym$, by means of $q$-commuting
variables~\cite{NTpfb}.

There is another way to obtain $\QSym$ from $\FQSym$: it is known~\cite{NCSF6}
that $\QSym$ is isomorphic to the image of $\FQSym(A)$ in the hypoplactic
algebra $\KK[A^*/\equiv_H]$.
One may then ask whether there exist a $q$-analogue of the hypoplactic
congruence leading directly to $\QSym_q$.

Recall that the hypoplactic congruence can be presented as the bi-sylvester
congruence:
\begin{equation}
\begin{split}
ubvcaw \equiv ubvacw, \text{ with } a<b\leq c \\
ucavbw \equiv uacvbw, \text{ with } a\leq b< c
\end{split}
\end{equation}
with $u,v,w\in A^*$.

A natural $q$-analogue, compatible with the above $q$-commutation is
\begin{equation}
\begin{split}
ubvcaw \equiv_{qH} q ubvacw, \text{ with } a<b\leq c \\
ucavbw \equiv_{qH} q uacvbw, \text{ with } a\leq b<c
\end{split}
\end{equation}
with $u,v,w\in A^*$.
Then, we have

\begin{theorem}
The image of $\FQSym(A)$ under the natural projection
$\KK\langle A\rangle \to \KK\langle A\rangle/\equiv_{qH}$ is isomorphic to
$\QSym_q$ as an algebra, and also as a twisted Hopf algebra for the coproduct
$A\to A'\dot+ A''$, $A'$ and $A''$ being $q$-commuting alphabets.
\end{theorem}

Moreover, if one only takes the sylvester congruence
\begin{equation}
ucavbw \equiv uacvbw,
\end{equation}
the quotient $\FQSym(A)$ under the natural projection
$\KK\langle A\rangle \to \KK\langle A\rangle/\equiv_S$ is isomorphic to the
Hopf algebra of planar binary trees of Loday and Ronco~\cite{LR1,HNT}. The
previous construction provides natural twisted $q$-analogues of this Hopf
algebra. Indeed, if one defines the $q$-sylvester congruence as
\begin{equation}
ucavbw \equiv_{qS} q uacvbw, \text{ with } a\leq b<c,
\end{equation}
then

\begin{theorem}
The image of $\FQSym(A)$ under the natural projection
$\KK\langle A\rangle \to \KK\langle A\rangle/\equiv_{qS}$ is a twisted Hopf
algebra, with basis indexed by planar binary trees.
\qed
\end{theorem}

\section*{Acknowledgements}
This project has been partially supported by EC's IHRP Programme, grant
HPRN-CT-2001-00272, ``Algebraic Combinatorics in Europe''.
The authors would also like to thank the contributors of the MuPAD project,
and especially of the combinat part, for providing the development environment
for this research.


\end{document}